\documentclass[pdftex]{article}
\usepackage{amsfonts, amsmath} 
\usepackage{xcolor}
\usepackage[colorlinks=true]{hyperref}
\hypersetup{linkcolor=blue}
\usepackage[dvips, lmargin=2.5cm, rmargin=1.5cm, tmargin=2.5cm, bmargin=2cm]{geometry}

\usepackage{tikz}
\tikzset{
	treenode/.style = {shape=, rounded corners,
		draw, align=center,
		top color=white,
		bottom color=blue!20},
	root/.style = {treenode, font=\Large,
		bottom color=red!30},
	env/.style = {treenode, font=\ttfamily\normalsize},
	dummy/.style = {circle,draw}
}
\usepackage{tikz}
\usetikzlibrary{arrows}
\tikzstyle{block}=[draw opacity=0.7,line width=1.4cm]
\newcommand{\verti}[1]{{\left\vert\kern-0.25ex\left\vert\kern-0.25ex\left\vert #1
    \right\vert\kern-0.25ex\right\vert\kern-0.25ex\right\vert}}
\title{
 \textbf{Stability analysis of SIR and SIRS models with non monotone incidence function and various mortality rates}
}

\numberwithin{equation}{section}

\newtheorem{theorem}{Theorem}

\newenvironment{prv}{ \noindent {\bf Proof.} \noindent} {\hfill$\square$ \vskip 6mm}
\date{}
\author{Y. Mohamed $^1$,  A. Ahmedou $^2$, M. S. B. Elemine Vall $^2$\\
$^1$ University of Nouakchott, Faculty of Legal and Economic Sciences, \\
Quantitative technics departement, Nouakchott, Mauritanie.\\
Email : {\color{blue}yahyajiddou@yahoo.fr}\\
 $^2$ University of Nouakchott, Professional university institute\\
Applied Mathematics and indistriel genious department,  Nouakchott, Mauritanie.\\
Email : {\color{blue}azizaahmedou@yahoo.fr}\\
Email : {\color{blue}saad2012bouh@gmail.com}
}
\everymath{\displaystyle}
\begin{document}
\maketitle
\vspace{0.5cm}
	\maketitle
	\begin{abstract}
		This study uses the Lyapunov method, the Poincar\'e-Bendixson theorem, and the Dulac criterion to analyze the stability of SIR and SIRS with non-monotone incidence and different mortality rates.
	\end{abstract}
	\section{Introduction}
	Mathematical models became important tools  for analyzing  the spread and control of infectious  diseases. The incidence function is an essential component to study  the spread of infectious disease. In general the incidence is represented as a linear function of the infection class, by a principle of mass action:
	\begin{equation}\label{blieair}
		g_1(I)S=\beta IS.
	\end{equation}
	see for instance  $\cite{An1, Bu1, Br1, Ke21, Ke22, Ke23}$ and the references therein.\\
	In their study of cholera epidemic, Capasso and Serio \cite{Ca1} were the first to use a nonlinear incidence of  Holling type II given by:
	\begin{equation}\label{HollingII}
		g_2(I)S=\frac{\beta IS}{1+\alpha I}
	\end{equation}
	where $\beta I$ measure the infection force of the disease and  $\frac{1}{1+\alpha I}$ measure the inhibition effect.\\
	In $\cite{Liu}$ Liu et al. was proposed a general incidence of the form:
	\begin{equation}\label{Gincidence}
		g(I)S=\frac{\beta I^pS}{1+\alpha I^q} 
	\end{equation}
	other authors have used this incidence see for example  $\cite{Heth}$ and $\cite{Der}.$\\
	If the the function $g(.)$ is non-monotone, that is  , $g(.)$ is increasing when $I$ is small and decreasing when $I$ is large, this incidence rate seems reasonable than the bi-linear incidence rate.
	
	We considered in this  paper a nonlinear and non-monotone incidence function of the previous form with $p=1$ and $q=2$, that is, 
	\begin{equation} \label{Inonmonotone}
		g(I)S=\frac{\beta IS}{1+\alpha I^2}
	\end{equation}
such incidence is used in \cite{Xi1}.\\
	The global stability is an interesting classic question in the mathematical epidemiology models. Korobeinikov and  Wake in $\cite{Ko23}$ used a Lyapunov function of Volterra type to show the global stability of SIS, SIR and SIRS models with bi-linear incidence.  Adda and Bichara have improved and completed Korobeinikov and  Wake results by considering some SIR and SIRS models with bi-linear incidence and differential mortality rates \cite{Ad2}. 
By using the Poincaré-Bendixson theorem, the Dulac criteria, and the Lyapunov function method to establish the conditions for global stability, Vargas-De-Leon has studied the global stability of the SIS, SIR, and SIRS models with standard incidence in \cite{art1, art2}. More specifically, the author has construct a Lyapunov function with two components, one of quadratic type and the second of Volterra type.

In this paper we propose two methods to show the global stability of SIR and SIRS models with nonmontone incidence rate and with various mortality rates. 
Our approach involves reducing the system to a planar one and then applying the Dulac criterion to eliminate periodic curves. To draw a conclusion, we also utilize the Poincaré-Bendixson theorem. Furthermore, we construct a Lyapunov function consisting of a quadratic term and a Volterra type term.
	
The paper is structured as follows: In the second section, we formulate the model and we calculate the basic reproduction number $R_0$ (It is defined as the number of new cases of infection caused by an infected individual in a population susceptible \cite{Di1} \cite{Va1}). In sections 3 and 4, we show that the disease-free equilibrium point is locally asymptotically stable if and only if $R_0<1$ and it is globally asymptotically stable. Finally when $R_0>1$ we show that the endemic equilibrium point is locally asymptotically stable and the global stability of endemic steady state using the Poincar\'e-Bendixson theorem, Dulac's criterion and Lyapunov's method in sections 5  and 6. In the section 7 we deduce the same previous results for a SIRS model. Finally, we conclude in section 8.
	\section{Model presentation}
	In this section we describe an SIR model with various mortality rate and non monotone incidence function by assuming that the incubation period is negligible.
	Let $N(t)$ be total population which is divided into three classes:
	susceptible individuals $S(t)$,  infected individuals $I(t)$, and individuals that have recovered/remove from the system $R(t)$. Based on this consideration, the total population is $N(t) = S(t) + I(t)+ R(t)$.
	
	We assume that all new born became susceptible and that the births recompense the deaths, that is,   
	\begin{equation}\label{eqmor}
	b= \mu_1S+\mu_2 I+\mu_3 R.
	\end{equation}
which means the size of the population $N(t)$ remains constant.

We model the contact by the following non monotone function:
	$$g(I)=\frac{\beta I}{1+\alpha I^2}$$
	with $\alpha >0$ , the terms $\beta I$ and $1+\alpha I^2$ represent respectively the force of the incidence and the inhibition effect. 
	
	We have the following scheme:
	
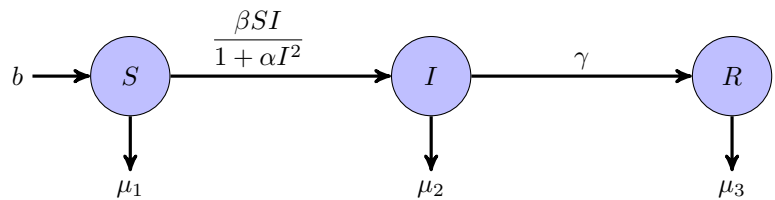
\begin{center}
		\begin{tikzpicture}
		\tikzset{radial/.style={very thick,->,>=stealth'}};
		\tikzstyle{point}=[circle,fill=blue!25,minimum width=3em]
		\tikzset{transversal/.style={->,>=stealth',thick,dashed}};
		\node[draw,point,rounded corners=3pt] (S)at(-2,0){$S$ };
		\node (SMu) at (-2,-1.5){$\mu_1$};
		\node (Sb) at (-3.5,0){$b$};
		\draw[radial] (S)--(SMu);
		\draw[radial] (Sb)--(S);
		\node[draw,point,rounded corners=3pt] (I)at(2,0) {$I$};
		\node[rounded corners=3pt] (F1)at(-0.3,0.5) {$\frac{\beta SI}{1+\alpha I^2}$};
		\draw[radial] (S)--(I);
		\node (IMu) at (2,-1.5){$\mu_2$};
		\draw[radial] (I)--(IMu);
		\node[draw,point,rounded corners=3pt] (R)at(6,0) {$R$ };
		\node[rounded corners=3pt] (F1)at(4,0.2) {$\gamma$};
		\draw[radial] (I)--(R);
		\node (RMu) at (6,-1.5){$\mu_3$};
		\draw[radial] (R)--(RMu);
		\end{tikzpicture}
	\end{center}
	
	Through the previous scheme, a system of nonlinear differential equations is obtained and presented below.
	\begin{equation} \label{a}
		\left\{\begin{array}{l}
			S^{\prime}=b-\mu_1 S-\frac{\beta S I}{(1+\alpha I^2)} . \\
			I^{\prime}=\frac{\beta S I}{(1+\alpha I^{2})}-(\mu_2+\gamma) I \\
			R^{\prime}=\gamma I-\mu_3 R .
		\end{array}\right.  
	\end{equation} 
	Thanks to (\ref{eqmor}) the previous system can be write
	\begin{equation}\label{b}
		\left\{\begin{array}{l}
			S^{\prime}=-\frac{\beta S I}{(1+\alpha I^2)}+ \mu_2 I +\mu_3 R . \\
			I^{\prime}=\frac{\beta S I}{(1+\alpha I^{2})}-(\mu_2+\gamma) I. \\
			R^{\prime}=\gamma I-\mu_3 R .
		\end{array}\right.  
	\end{equation}
	Since  $S(t)+I(t)+R(t)=N $=Constant, the system $(\ref{b})$ is equivalent to the following  planar system :
	
	\begin{equation}\label{c}
		\left\{\begin{array}{l}
			S^{\prime}=-\frac{\beta S I}{(1+\alpha I^2)}+ \mu_2 I +\mu_3 (N-S-I) . \\
			I^{\prime}=\frac{\beta S I}{(1+\alpha I^{2})}-(\mu_2+\gamma) I. \\
		\end{array}\right. 
	\end{equation}
	
	For the simplicity of notations, we put   $S=\frac{S}{N}$, $I=\frac{I}{N}$, $\beta=N\beta$ and $\alpha=N^2\alpha$, which gives the following system 
	\begin{equation} \label{dd}
		\left\{\begin{array}{l}
			S^{\prime}= \mu_3 +(\mu_2 -\mu_3)I- \mu_3 S -\frac{\beta S I}{1+\alpha I^2} . \\
			I^{\prime}=\frac{\beta S I}{1+\alpha I^{2}}-(\mu_2+\gamma) I. \\
		\end{array}\right. 
	\end{equation}
	with  $S\geq 0, I \geq 0$  such that  $ S+I \leq 1.$ 
	
	The biological domain of the previous system is the following standard compact positively invariant simplex given by 
	$$D= \left\{  (S,I)\in \mathbb{R}^2_{+}: S+I \leq 1\right\}$$
	\section{Local stability of the free-diseases equilibrium}
	It is very easy to show that the system $(\ref{dd})$ admits a single free-diseases equilibrium point given by $x_{DFE}=(S_0,0)=(1,0).$\\
	Furthermore, the basic reproduction number of the system $(\ref{dd})$ is obtained by  Diekmann et al. in \cite{Di1} and Driessche and Watmough in \cite{Va1} as follows
	$$ R_0= \frac{1}{\mu_2 + \gamma} \times \frac{\partial g}{\partial I}(1,0)=\frac{\beta}{\mu_2 + \gamma}.$$
	This enables us to state the following local stability result.
	\begin{theorem}\label{th1}
		The free-diseases equilibrium point $x_{DFE}$ is locally asymptotically stable if we assume that $R_0 <1$ and unstable elsewhere.
	\end{theorem}
	\begin{prv}
		This proof based on the direct calculation of eigenvalues of the Jacobin matrix of system, that is, the matrix defined by
		\begin{eqnarray*}
			J(S,I)=
			\begin{bmatrix}
				- \mu_3-\frac{\beta I}{1+\alpha I^2} & (\mu_2-\mu_3)-\frac{\beta S(1-\alpha I^2)}{(1+\alpha I^2)^2}  \\
				\frac{\beta I}{1+\alpha I^2} & \frac{\beta S(1-\alpha I^2)}{(1+\alpha I^2)^2}-(\mu_2 +\gamma)  \\
				
			\end{bmatrix}
		\end{eqnarray*}
		at the free-diseases equilibrium point, we have 
		\begin{eqnarray*}
			J(1,0)=
			\begin{bmatrix}
				- \mu_3 & \mu_2-\mu_3-\beta  \\
				0 & \beta-(\mu_2 +\gamma)  \\
			\end{bmatrix}
		\end{eqnarray*}
		Therefore the eigenvalues of this matrix are  $\lambda_1= -\mu_3 <0$ and $\lambda_2=(\mu_2 +\gamma)(R_0-1)$.\\
		Then, if $R_0\geq 1$   \ the second eigenvalue $\lambda_2$ remains positive, which means that  the free-diseases equilibrium point $x_{DFE}$ is unstable and if  $R_0 <1$, the second eigenvalue $\lambda_2<0$ then $x_{DFE}$ is asymptotically stable. 
		
	\end{prv}
	
	\section{Global stability of free-diseases equilibrium point}
	\begin{theorem}
		Under the hypothesis $R_{0}<1$, the free-diseases equilibrium point $x_{DFE}$  is globally asymptotically stable.
	\end{theorem}
	\begin{prv}
		Consider the following Lyapunov function
		
		\begin{eqnarray*}
			V(S,I) =  I.
		\end{eqnarray*}
		Then \begin{eqnarray*}
			V'(S,I) &= & I'\\
			&=&\frac{\beta S I}{1+\alpha I^{2}}-(\mu_2+\gamma) I\\
			& \leq & (\mu_2+\gamma)(R_0 S-1)I\\
			&\leq & 0
		\end{eqnarray*}
		If we assume that $V'=0$ then $I=0$ or ($S=S_0$ and $R_0=1.$)  So the largest invariant set contained in the biological domain $D$ is $M=\left\{ (S,I)\in D: V'=0\right\}$ which is reduced to the free-diseases equilibrium point $x_{DFE}$ and by the LaSalle  invariance principle in  $\cite{La21, La22}$ the point $x_{DFE}$ is globally asymptotically stable.
	\end{prv}
	\section{Local stability of the endemic equilibrium point}
	We look at the existence of endemic equilibrium point. This endemic equilibrium state satisfies
	$$
	\left\{\begin{array}{l}
	0= \mu_3 +(\mu_2 -\mu_3)I- \mu_3 S^\ast -\frac{\beta S^\ast I^\ast}{1+\alpha (I^\ast)^2} . \\
	0=\frac{\beta S^\ast I^\ast}{1+\alpha (I^\ast)^{2}}-(\mu_2+\gamma) I^\ast. \\
	\end{array}\right.
	$$
	By solving the previous nonlinear system, we will have
	$$x^{*}=(S^\ast,I^\ast)=\left(\frac{1+\alpha{I^\ast}^2}{R_0}, \quad \frac{-(\mu_3+\gamma)R_0+\sqrt{(\mu_3+\gamma)^2R_0^2+4\alpha\mu_3^2(R_0-1)}}{2\alpha\mu_3}\right).$$
	It is easy to  show that if $R_0>1,$ then $S^{*}>0, I^{*}>0$ and $S^{*}+I^{*} \leq 1$.	
	Now, we are in position to quote the following local stability theorem.
	\begin{theorem}\label{las}
		If $R_{0}>1$. Then, the endemic equilibrium point  $x^{*}$ is locally asymptotically stable.
	\end{theorem}
	\begin{prv} 
		This proof is quite similar to the one in theorem \ref{th1}. We start by giving the Jacobin matrix of system (\ref{dd}) at the endemic equilibrium, that is ,
		\begin{eqnarray*}
			J(S^\ast,I^\ast)=
			\begin{bmatrix}
				- \mu_3-\frac{\beta I^*}{1+\alpha {I^\ast}^2} & (\mu_2-\mu_3)-\frac{\beta S^*(1-\alpha {I^\ast}^2)}{(1+\alpha {I^\ast}^2)^2}  \\
				\frac{\beta I^*}{1+\alpha {I^\ast}^2} & \frac{\beta S^*(1-\alpha {I^\ast}^2)}{(1+\alpha {I^\ast}^2)^2}-(\mu_2 +\gamma)  \\
				
			\end{bmatrix}
		\end{eqnarray*}
		Since at the endemic equilibrium point, we have 
		$$\frac{\beta S^{*} I^{*}}{1+\alpha {I^*}^{2}}=(\mu_2+\gamma) I^{*}$$
		which yields $$\mu_2+\gamma=\frac{\beta S^{*}}{1+\alpha {I^*}^{2}}$$
		and $$\frac{\beta S^*-\alpha \beta S^*  {I^*}^{2}}{\left(1+\alpha {I^*}^{2}\right)^{2}}-(\mu_2+\gamma)=\frac{ \beta S^*-\alpha \beta S^*{I^*}^{2}}{\left(1+\alpha {I^*}^{2}\right)^{2}}-\frac{\beta S^*\left(1+\alpha {I^*}^{2}\right)}{\left(1+\alpha {I^*}^{2}\right)^{2}}=\frac{-2\alpha \beta S^{*} {I^*}^{2}}{\left(1+\alpha {I^*}^{2}\right)^{2}},$$
		and because 
		$$ \mu_2-\mu_3= \frac{\beta S^*}{1+\alpha {I^\ast}}-\frac{\mu_3}{I^*}(1-S^*)$$
		The Jacobain matrix $J(x^*)$, can be write
		\begin{eqnarray*}
			J(x^*)=
			\begin{bmatrix}
				- \mu_3-X & Y-\frac{\mu_3}{I^*}(1-S^*) \\
				X & -Y \\
			\end{bmatrix}
		\end{eqnarray*}
		where $X=\frac{\beta I^*}{1+\alpha {I^\ast}^2},$ $Y=\frac{\alpha \beta S^{*} {I^*}^{2}}{2\left(1+\alpha {I^*}^{2}\right)^{2}}$
		and in consequences, we have
		$$Tr(J(x^*))= -\mu_3-X-Y<0\quad\textrm{ and }\quad\det J(x^*)= \mu_3 Y +\frac{\mu_3}{I^*}(1-S^*)X \geq 0$$
		Thus, the eigenvalues of $(J(x^*))$ have a negative real part, which means the endemic equilibrium $x^*$ is locally asymptotically stable.
	\end{prv}
	\section{Global stability of endemic equilibrium point}
	In this section, we propose two methods to prove the global stability result state in the following theorem.
	\subsection{Dulac's criteria}
	\begin{theorem}
		Assume that $R_0>1.$ Then, the endemic equilibrium point is globally asymptotically stable. 
	\end{theorem}
	\begin{prv}
		Consider the following planar system
		$$\left\{\begin{array}{l}
		S^{\prime}= \mu_3 +(\mu_2 -\mu_3)I- \mu_3 S -\frac{\beta S I}{1+\alpha I^2}=f_1(S,I) . \\
		I^{\prime}=\frac{\beta S I}{1+\alpha I^{2}}-(\mu_2+\gamma) I=f_2(S,I). \\
		\end{array}\right.  \label{d}$$ 
		We employ the Dulac's criterion to achieve this proof by taking the following Dulac function $B(S,I)=\frac{1+\alpha I^{2}}{\beta SI}$, we obtain
		$$\frac{\partial Bf_1}{\partial S}(S,I)+\frac{\partial Bf_2}{\partial I}(S,I)=-\frac{\mu_2 (1+\alpha I^{2}) }{\beta S^2}-\frac{2\alpha I(\mu_2+\gamma)}{\beta S}+\frac{\mu_3(1+\alpha I^{2})}{\beta S^2}(1-\frac{1}{I})<0$$
		Thus, the system (\ref{dd}) does not have a limit cycle in $D$.
		
		In view of theorem \ref{las}  $x^*$ is locally asymptotically stable and since $D$ is positively invariant set, then Poincar\'e-Bendixson theorem shows that   $x^*$ is globally asymptotically stable.
	\end{prv}
	\subsection{Lyapunov method}	
	Now, we are in position to state another proof of theorem \ref{th1} based on Lyapunov functions.
	
	\begin{prv}
		Let $V$ be the Lyapunov function defined by
		$$ V(S,I)=\frac{1}{2}(S-S^\ast+I-I^\ast)^2+\frac{2\mu_3+\gamma}{\beta}(1+\alpha {I^\ast}^2)\left(I-I^\ast-I^\ast\ln\left(\frac{I}{I^\ast}\right)\right).$$
		Then, the time derivative of $V$ may be write as
		\begin{eqnarray*}
			V'(S,I)&=&(I'+S')(S-S^\ast+I-I^\ast)+\frac{2\mu_3+\gamma}{\beta}(1+\alpha {I^\ast}^2)I'\left(1-\frac{I^\ast}{I}\right)\\
			&=&(\mu_3-\mu_3S-(\mu_3+\gamma)I)(S-S^\ast+I-I^\ast)+\frac{2\mu_3+\gamma}{\beta}(1+\alpha {I^\ast}^2)\left(\frac{\beta SI}{1+\alpha I^2}-(\mu_2+\gamma)I)\right)\left(1-\frac{I^\ast}{I}\right)\\
			&=&(\mu_3(S^\ast-S)-(\mu_3+\gamma)(I^\ast-I))(S-S^\ast+I-I^\ast)+(2\mu_3+\gamma)(1+\alpha {I^\ast}^2)\left(\frac{S}{1+\alpha I^2}-\frac{S^\ast}{1+\alpha {I^\ast}^2}\right)\left(I-I^\ast\right)\\
			&=&-\mu_3(S-S^\ast)^2-(\mu_3+\gamma)(I-I^\ast)^2-(2\mu3+\gamma)(I-I^\ast)(S-S^\ast)\\
			&&+(2\mu_3+\gamma)(1+\alpha{I^\ast}^2)\left[\frac{S-S^\ast}{1+\alpha {I^\ast}^2}-\frac{\alpha S({I^\ast}^2-I^2)}{(1+\alpha I^2)(1+\alpha {I^\ast}^2)}\right](I-I^\ast)\\
			&=&-\mu_3(S-S^\ast)^2-(\mu_3+\gamma)(I-I^\ast)^2-\frac{(2\mu_3+\gamma)\alpha S(I^\ast+I)}{1+\alpha I^2}(I-I^\ast)^2\\
			&\leq&0.
		\end{eqnarray*} 
	
Therefore, $V'=0$ if and only if $S=S^\ast$ and $I=I^\ast$. Thus the largest set $\{(S,I)\in D: V'=0\}$ 	is reduced to $x^\ast$.
			
Then, according to LaSalle principle in $\cite{La21, La22}$ the positive equilibrium state $x^\ast$ is globally asymptotically stable in $D.$
	\end{prv}
	\section{SIRS model}
	In this section, we consider a SIRS  with non monotone incidence and with different mortality. We keep the same notations in the sections above, we have  the following system:
	\begin{equation}\label{SIRS}
		\left\{\begin{array}{l}
			S^{\prime}=b-\mu_1 S-\frac{\beta S I}{(1+\alpha I^2)}+\rho R . \\
			I^{\prime}=\frac{\beta S I}{(1+\alpha I^{2})}-(\mu_2+\gamma) I \\
			R^{\prime}=\gamma I-(\mu_3 +\rho) R .
		\end{array}\right. 
	\end{equation}
	which reduces to 
	\begin{equation}\label{reduitSIRS}
		\left\{\begin{array}{l} 
			S^{\prime}=-\frac{\beta S I}{(1+\alpha I^2)}+ \mu_2 I +(\mu_3 +\rho)R . \\
			I^{\prime}=\frac{\beta S I}{(1+\alpha I^{2})}-(\mu_2+\gamma) I. \\
			R^{\prime}=\gamma I-(\mu_3 +\rho)R .
		\end{array}\right.  
	\end{equation}
	The system $(\ref{reduitSIRS})$ is exactly as system $(\ref{b})$ where $\mu_3$ replaced by $\mu_3+\rho.$
	\section{Conclusion}
	In this study, we established the global stability of SIR and SIRS epidemiological models with non monotone incidence and varied mortality rates by using the Poincar\'e-Bendixson theorem, Dulac's criterion, and Lyapunov approach. The results of the study from Adda and Bichara in \cite{Ad2} and Korobeinikov and Wake in \cite{Ko23} were generalized in this work.
	
\end{document}